\documentclass[a4paper, 12pt]{article}
\pdfoutput=1 
\usepackage[T1]{fontenc}
\usepackage{lmodern}
\usepackage{etoolbox}
\usepackage{amsmath, amsthm, amssymb, bm, mleftright}
\usepackage{booktabs, longtable, caption, subcaption, multirow}
\usepackage[space]{grffile}
\usepackage[margin=2.5cm, footskip=1cm]{geometry}
\usepackage{enumitem}
\setenumerate[1]{label=(\arabic*), ref=\arabic*}
\usepackage[protrusion=allmath]{microtype}
\usepackage{tikz}
\usetikzlibrary{cd, calc}
\usepackage{listings}
\usepackage[colorlinks, citecolor=blue, linkcolor=blue, linktoc=section]{hyperref}

\allowdisplaybreaks
\usepackage[capitalise, noabbrev]{cleveref}
\creflabelformat{enumi}{(#2#1#3)}
\creflabelformat{enumii}{(#2#1#3)}
\creflabelformat{enumiii}{(#2#1#3)}

\usepackage{titlesec}
\titleformat*{\section}{\large\bfseries}
\titleformat*{\subsection}{\normalsize\bfseries}
\newlength{\VerticalSpaceAfterParagraph}
\titlespacing*{\paragraph}{0pt}{\VerticalSpaceAfterParagraph}{1em}

\usepackage{titling}
\pretitle{\vspace{-\baselineskip}\begin{center}\Large\bfseries}
\posttitle{\end{center}\vspace{-0.25\baselineskip}}
\preauthor{\begin{center}}
\postauthor{\end{center}}
\predate{\begin{center}}
\postdate{\end{center}}

\setlist
  {
    topsep = 5.0pt plus 2.0pt minus 3.0pt,
    partopsep = 1.5pt plus 1.0pt minus 1.0pt,
    parsep = 2.5pt plus 1.25pt minus 0.5pt,
    itemsep = 0pt plus 1.25pt minus 0.5pt
  }

\theoremstyle{definition}
\newtheorem{theorem}{Theorem}
\newtheorem{proposition}[theorem]{Proposition}
\newtheorem{question}[theorem]{Question}
\newtheorem{definition}[theorem]{Definition}
\newtheorem{problem}[theorem]{Problem}

\setenumerate[1]{label=\textup{\arabic{enumi}.}, ref=\arabic{enumi}}
\setenumerate[2]{label=\textup{(\alph*)}, ref=\alph*}

\creflabelformat{enumi}{#2#1#3}
\creflabelformat{enumii}{#2#1#3}

\Crefname{problem}{Problem}{Problems}

\numberwithin{theorem}{section}
\numberwithin{table}{section}
\numberwithin{equation}{theorem}

\author{Erik Paemurru}

\title{Extended Abstract of Counting Divisorial Contractions with Centre a $cA_n$-Singularity}
\date{\vspace{-2\baselineskip}}

\begin{document}

\maketitle

\begin{abstract}
This text is chapter 149 of the ZAG Handbook of Algebraic Geometry. We review the main results of the paper ``Counting divisorial contractions with centre a $cA_n$-singularity'', published in Publications of the Research Institute for Mathematical Sciences, Kyoto University.
\end{abstract}

\setcounter{section}{149}

\subsection{Introduction}

A \emph{variety} is an integral separated scheme of finite type over the complex numbers~$\mathbb C$.

\begin{definition}
A \emph{divisorial contraction} is a proper birational morphism between normal varieties (or a proper bimeromorphic holomorphic map between normal integral complex spaces) with terminal singularities such that the exceptional locus is a relatively antiample $\mathbb Q$-Cartier prime divisor.
\end{definition}

When the base variety is $\mathbb Q$-factorial and proper over~$\mathbb C$, then the divisorial contraction is necessarily the contraction of an extremal ray. Divisorial contractions appear naturally in the minimal model program and the Sarkisov program. We do not require $\mathbb Q$-factoriality in the definition because three-dimensional divisorial contractions have been studied locally analytically and $\mathbb Q$-factoriality is not always preserved when restricting to analytic open subsets.

The explicit classification of three-dimensional divisorial contractions has been used in results on non-rationality, birational rigidity, birational non-rigidity and birational solidity. See for example \cite{CM04,Oka23,KOPP24,AK16,ACP21}. In the mentioned papers, the divisorial contraction is either used to produce a Sarkisov link or, after a computation involving intersection numbers, it is shown that the divisorial contraction cannot initiate a Sarkisov link.

\subsection{Ultimate description}

See \cite[Definition~2.3.2 and Definition~2.2.10]{Kaw24} for the definition of compound $A_k$ singularities, denoted $cA_k$, and quotient singularities. By \cite[Theorem~0.6]{Rei83}, Gorenstein three-dimensional terminal singularities are isolated \emph{compound Du Val singularities}.

We say two proper birational morphisms are \emph{locally analytically equivalent} if there exists a local biholomorphism around the centre that lifts to a local biholomorphism around the exceptional loci on the blown-up spaces (\cite[Definition~3.3]{Pae24Counting}).
We have the following result, which is called the \emph{ultimate description} in \cite[Theorem~3.5.8]{Kaw24}:
\begin{theorem} \label{thm:ultimate description}
Let $\varphi$ be a three-dimensional divisorial contraction which contracts a prime divisor to a point $P$ such that one of the following holds:
\begin{enumerate}
\item the centre $P$ of $\varphi$ is a non-Gorenstein point, a smooth point or a $cA_k$ point, or
\item the discrepancy of the exceptional divisor of $\varphi$ is greater than~$1$.
\end{enumerate}
Then $\varphi$ is locally analytically equivalent to a suitable weighted blowup of
\[
V\bigl(f(x_1, x_2, x_3, x_4, x_5),~ g(x_1, x_2, x_3, x_4) + x_5\bigr) \subseteq \frac{1}{n}(a_1, a_2, a_3, a_4, a_5)
\]
with respect to the orbifold coordinates $x_1, x_2, x_3, x_4, x_5$.
\end{theorem}
Moreover, the weights and the convergent power series $f$ and $g$ in \cref{thm:ultimate description} are described explicitly. The case of non-Gorenstein point in \cref{thm:ultimate description} is due to \cite{Hay99,Hay00,Hay05,Kaw05,Kaw12,Kaw96} and the other cases are due to \cite{Kaw01,Kaw02,Kaw03}, \cite[Theorem~1.2]{Kaw05} and \cite{Yam18}.

\subsection{Simplifying the classification}

The above papers of Kawamata, Kawakita, Hayakawa and Yamamoto give a list of morphisms such that every three-dimensional divisorial contraction, except for the aforementioned $cD_k$ and $cE_k$ case, is locally analytically equivalent to one of the morphisms in the list. For $cA_k$ singularities, we have the following \lcnamecref{thm:classification for cAk}.

\begin{theorem}[{\cite[Theorem~1.1]{Kaw02}, \cite[Theorem~1.13]{Kaw03} and \cite[Theorem~6.1]{Pae24Counting}}] \label{thm:classification for cAk}
Let $k$ be a positive integer and $P$ a $cA_k$ singularity of a terminal variety~$X$. Let $\varphi\colon Y \to X$ be a proper birational morphism with centre~$P$. Then $\varphi$ is a divisorial contraction if and only if one of the following holds:
\begin{enumerate}
\item \label{itm:pre div contr cAn - gen} $\varphi$ is locally analytically equivalent to the $(r_1, r_2, a, 1)$-blowup of $\mathbb V(f)$ at $\bm 0$ where $f \in \mathbb C\{x, y, z, t\}$ is such that
  \begin{enumerate}
  \item \label{itm:integers} $r_1, r_2$ and $a$ are positive integers such that $r_1 \leq r_2$, $a(k+1) = r_1 + r_2$ and $a$ is coprime to both $r_1$ and $r_2$, and
  \item \label{itm:pre div contr cAn - gen f} $f = xy + g(z, t)$ where $\operatorname{wt} g = r_1 + r_2$,
  \end{enumerate}
\item \label{itm:pre div contr cAn - cA1} $k = 1$ and $\varphi$ is locally analytically equivalent to the $(1, 5, 3, 2)$-blowup of $\mathbb V(f)$ at $\bm 0$ where $f \in \mathbb C\{x, y, z, t\}$ is such that
  \begin{enumerate}
  \item \label{itm:pre div contr cAn - cA1 f} $f = xy + z^2 + t^3$,
  \end{enumerate}
\item \label{itm:pre div contr cAn - cA2} $k = 2$ and $\varphi$ is locally analytically equivalent to the $(4, 3, 2, 1)$-blowup of $\mathbb V(f)$ at $\bm 0$ where $f \in \mathbb C\{x, y, z, t\}$ is such that
  \begin{enumerate}
  \item \label{itm:pre div contr cAn - cA2 f} $f = x^2 + y^2 + z^3 + x t^2$.
  \end{enumerate}
\end{enumerate}
\end{theorem}

Item \labelcref{itm:pre div contr cAn - cA2} is the simplification \cite[Theorem~6.1(3)]{Pae24Counting} of \cite[Theorem~2.6]{Yam18}, which had the following complicated condition instead of \labelcref{itm:pre div contr cAn - cA2}\labelcref{itm:pre div contr cAn - cA2 f}:
\begin{enumerate}
\item[($*$)] $f = x^2 + y^2 + 2cxy + 2xp(z, t) + 2cyp_{\operatorname{wt}=3}(z, t) + z^3 + g(z, t)$, where $c \in \mathbb C \setminus \{-1, 1\}$, $\operatorname{wt} g \geq 6$, the power series $p$ contains only monomials of weight $2$ and $3$, the coefficient of $t^2$ is non-zero in $p$ and $\deg g(z, 1) \leq 2$.
\end{enumerate}

Ideally, we would like a stronger classification:

\begin{problem}[{\cite[Problem~1.1]{Pae24Counting}}] \label{pro:int loc ana equiv classes}
Describe the local analytic equivalence classes of three-dimensional divisorial contractions with centre a point.
\end{problem}

For $cA_k$ singularities, this follows from the \lcnamecref{thm:local analytic equivalence classes} below.

\begin{proposition}[{\cite[Proposition~4.7]{Pae24Sextic}}] \label{thm:local analytic equivalence classes}
Let $k$ be a positive integer and let the variables $x, y, z, t$ have weights $(r_1, r_2, a, 1)$ as in \cref{thm:classification for cAk} part \labelcref{itm:pre div contr cAn - gen}\labelcref{itm:integers}.
Let $f, f' \in \mathbb C\{x, y, z, t\}$ define locally biholomorphic $cA_k$ singularities at the origin.
If $\operatorname{wt} f = \operatorname{wt} f' = r_1 + r_2$, then there exists a biholomorphic map germ $(V(f), \bm 0) \to (V(f'), \bm 0)$ that lifts to a local biholomorphism around the exceptional divisors of the weighted blown-up spaces.
\end{proposition}

Given a divisorial contraction, we would like to be able to determine which local analytic equivalence class it belongs to, ideally by comparing a finite set of discrete invariants.

\begin{problem}[{\cite[Problem~1.3]{Pae24Counting}}] \label{pro:int algorithm}
Describe an algorithm to determine whether a given weighted blowup is locally analytically equivalent to a given member of the classification list.
\end{problem}

For $cA_k$ singularities, this is provided by the next \lcnamecref{thm:swap}.

\begin{theorem}[{\cite[Theorem~6.1]{Pae24Counting}}] \label{thm:swap}
\Cref{thm:classification for cAk} continues to hold when we do all of the following substitutions:
\begin{itemize}
\item swap \labelcref{itm:pre div contr cAn - gen}\labelcref{itm:pre div contr cAn - gen f} with ``$\operatorname{wt} f = r_1 + r_2$'',
\item swap \labelcref{itm:pre div contr cAn - cA1}\labelcref{itm:pre div contr cAn - cA1 f} with ``$(V(f), \bm 0)$ is an $A_2$ singularity and $\operatorname{wt} f = 6$'', and
\item swap \labelcref{itm:pre div contr cAn - cA2}\labelcref{itm:pre div contr cAn - cA2 f} with ``$(V(f), \bm 0)$ is an $E_6$ singularity and $\operatorname{wt} f = 6$''.
\end{itemize}
\end{theorem}

\subsection{Global algebraic classification}

To construct Sarkisov links, we would like a global algebraic description of divisorial contractions. Assuming we have a local analytic classification, the global algebraic classification is provided by the following:

\begin{proposition}[{\cite[Corollary~5.6]{Pae24Counting}}] \label{thm:global algebraic}
Let $U = \operatorname{Spec} (\mathbb C[x_1, \ldots, x_n] / I)$ be an affine variety properly containing the origin. Assign positive integer weights $w_1, \ldots, w_m$ to the variables $y_1, \ldots, y_m$ of $\mathbb C^m$ and assign the weights $1, \ldots, 1$ to $x_1, \ldots, x_n$. Let superscript ``an'' denote analytification. Let $\psi\colon (U^{\mathrm{an}}, \bm 0) \to (Z, \bm 0)$ be a biholomorphic map germ to a complex space subgerm $(Z, \bm 0)$ of $(\mathbb C^m, \bm 0)$. Define the variety $\hat U$ by
\[
\hat U := \operatorname{Spec} \frac{\mathbb C[x_1, \ldots, x_n, y_1, \ldots, y_m]}{I + (\psi_1^{<w_1} - y_1, \ldots, \psi_m^{<w_m} - y_m)},
\]
where $\psi_j^{<w_j}$ denotes the truncation of the $j$-th coordinate power series of $\psi$ up to order $w_j - 1$. Note that $\hat U$ is isomorphic to~$U$ and let $\theta\colon (\hat U^{\mathrm{an}}, \bm 0) \to (Z, \bm 0)$ be the composition of the map germs of $\hat U^{\mathrm{an}} \to U^\mathrm{an}$ and $U^{\mathrm{an}} \to Z$. Then, both of the following hold:
\begin{enumerate}
\item $\theta$ lifts to a biholomorphism around the exceptional loci of the weighted blowups, and
\item every proper birational morphism locally analytically equivalent to the weighted blowup of $Z$ is given by the weighted blowup of $\hat U$ for some local biholomorphism~$\psi\colon U \to Z$.
\end{enumerate}
\end{proposition}

\subsection{Counting divisorial contractions}

For the explicit construction of all divisorial contractions, for example for the purpose of constructing all Sarkisov links, it is helpful to know whether there exist finitely or uncountably many divisorial contractions.

\begin{question}[{\cite[Question~1.5]{Pae24Counting}}]
Let $X$ be a three-dimensional terminal variety and $P \in X$ a singular point. Do there exist only finitely many divisorial contractions to $X$ with centre $P$ up to local analytic equivalence?
\end{question}

By \cite[Theorem~6.5]{Pae24Counting}, the answer is ``yes'' if $P$ is a $cA_k$ singularity.

For a given variety $X$ and $cA_k$ singularity $P \in X$, if there exists one divisorial contraction to $X$ with centre $P$ with a discrepancy of at least two, then by \cite[Theorem~6.5(c)]{Pae24Counting}, globally algebraically there exist uncountably many divisorial contractions to $X$ with centre~$P$. This statement does not always hold if $P$ is a $cE_7$ point, as shown in \cite[Remark~6.17]{Oka23}. This answers \cite[Question~1.6]{Pae24Counting} negatively.

\providecommand{\bysame}{\leavevmode\hbox to3em{\hrulefill}\thinspace}
\providecommand{\MR}{\relax\ifhmode\unskip\space\fi MR }
\providecommand{\MRhref}[2]{%
  \href{http://www.ams.org/mathscinet-getitem?mr=#1}{#2}
}
\providecommand{\href}[2]{#2}

\end{document}